\documentclass[12pt]{amsart}
\usepackage{amsmath,amssymb,amsbsy,amsfonts,amsthm,latexsym,amsopn,amstext,
                                                       amsxtra,euscript,amscd}
\begin{document}

\def\bbbr{{\r\m I\!R}}
\def\bbbc{{\rm I\!C}}
\def\bbbm{{\rm I\!M}}
\def\bbbn{{\rm I\!N}}
\def\bbbf{{\rm I\!F}}
\def\bbbh{{\rm I\!H}}
\def\bbbk{{\rm I\!K}}
\def\bbbl{{\rm I\!L}}
\def\bbbp{{\rm I}}
\def\sssum{\mathop{\sum\!\sum\!\sum}}
\newcommand{\lcm}{{\rm lcm}}
\newcommand{\eq}{{\mathbf{e}_q}}
\newcommand{\Em}{{\mathbf{e}_m}}
\def\bbbone{{\mathchoice {\rm 1\mskip-4mu l} {\rm 1\mskip-4mu l}
{\rm 1\mskip-4.5mu l} {\rm 1\mskip-5mu l}}}
\def\bbbc{{\mathchoice {\setbox0=\hbox{$\displaystyle\rm C$}\hbox{\hbox
to0pt{\kern0.4\wd0\vrule height0.9\ht0\hss}\box0}}
{\setbox0=\hbox{$\textstyle\rm C$}\hbox{\hbox
to0pt{\kern0.4\wd0\vrule height0.9\ht0\hss}\box0}}
{\setbox0=\hbox{$\scriptstyle\rm C$}\hbox{\hbox
to0pt{\kern0.4\wd0\vrule height0.9\ht0\hss}\box0}}
{\setbox0=\hbox{$\scriptscriptstyle\rm C$}\hbox{\hbox
to0pt{\kern0.4\wd0\vrule height0.9\ht0\hss}\box0}}}}
\def\bbbq{{\mathchoice {\setbox0=\hbox{$\displaystyle\rm
Q$}\hbox{\raise
0.15\ht0\hbox to0pt{\kern0.4\wd0\vrule height0.8\ht0\hss}\box0}}
{\setbox0=\hbox{$\textstyle\rm Q$}\hbox{\raise
0.15\ht0\hbox to0pt{\kern0.4\wd0\vrule height0.8\ht0\hss}\box0}}
{\setbox0=\hbox{$\scriptstyle\rm Q$}\hbox{\raise
0.15\ht0\hbox to0pt{\kern0.4\wd0\vrule height0.7\ht0\hss}\box0}}
{\setbox0=\hbox{$\scriptscriptstyle\rm Q$}\hbox{\raise
0.15\ht0\hbox to0pt{\kern0.4\wd0\vrule height0.7\ht0\hss}\box0}}}}
\def\bbbt{{\mathchoice {\setbox0=\hbox{$\displaystyle\rm
T$}\hbox{\hbox to0pt{\kern0.3\wd0\vrule height0.9\ht0\hss}\box0}}
{\setbox0=\hbox{$\textstyle\rm T$}\hbox{\hbox
to0pt{\kern0.3\wd0\vrule height0.9\ht0\hss}\box0}}
{\setbox0=\hbox{$\scriptstyle\rm T$}\hbox{\hbox
to0pt{\kern0.3\wd0\vrule height0.9\ht0\hss}\box0}}
{\setbox0=\hbox{$\scriptscriptstyle\rm T$}\hbox{\hbox
to0pt{\kern0.3\wd0\vrule height0.9\ht0\hss}\box0}}}}
\def\bbbs{{\mathchoice
{\setbox0=\hbox{$\displaystyle     \rm S$}\hbox{\raise0.5\ht0\hbox
to0pt{\kern0.35\wd0\vrule height0.45\ht0\hss}\hbox
to0pt{\kern0.55\wd0\vrule height0.5\ht0\hss}\box0}}
{\setbox0=\hbox{$\textstyle        \rm S$}\hbox{\raise0.5\ht0\hbox
to0pt{\kern0.35\wd0\vrule height0.45\ht0\hss}\hbox
to0pt{\kern0.55\wd0\vrule height0.5\ht0\hss}\box0}}
{\setbox0=\hbox{$\scriptstyle      \rm S$}\hbox{\raise0.5\ht0\hbox
to0pt{\kern0.35\wd0\vrule height0.45\ht0\hss}\raise0.05\ht0\hbox
to0pt{\kern0.5\wd0\vrule height0.45\ht0\hss}\box0}}
{\setbox0=\hbox{$\scriptscriptstyle\rm S$}\hbox{\raise0.5\ht0\hbox
to0pt{\kern0.4\wd0\vrule height0.45\ht0\hss}\raise0.05\ht0\hbox
to0pt{\kern0.55\wd0\vrule height0.45\ht0\hss}\box0}}}}
\def\bbbz{{\mathchoice {\hbox{$\sf\textstyle Z\kern-0.4em Z$}}
{\hbox{$\sf\textstyle Z\kern-0.4em Z$}}
{\hbox{$\sf\scriptstyle Z\kern-0.3em Z$}}
{\hbox{$\sf\scriptscriptstyle Z\kern-0.2em Z$}}}}

\newtheorem{theorem}{Theorem}
\newtheorem{lemma}[theorem]{Lemma}
\newtheorem{claim}[theorem]{Claim}
\newtheorem{cor}[theorem]{Corollary}
\newtheorem{prop}[theorem]{Proposition}
\newtheorem{definition}{Definition}
\newtheorem{question}[theorem]{Open Question}

\def\squareforqed{\hbox{\rlap{$\sqcap$}$\sqcup$}}
\def\qed{\ifmmode\squareforqed\else{\unskip\nobreak\hfil
\penalty50\hskip1em\null\nobreak\hfil\squareforqed
\parfillskip=0pt\finalhyphendemerits=0\endgraf}\fi}

\def\cA{{\mathcal A}}
\def\cB{{\mathcal B}}
\def\cC{{\mathcal C}}
\def\cD{{\mathcal D}}
\def\cE{{\mathcal E}}
\def\cF{{\mathcal F}}
\def\cG{{\mathcal G}}
\def\cH{{\mathcal H}}
\def\cI{{\mathcal I}}
\def\cJ{{\mathcal J}}
\def\cK{{\mathcal K}}
\def\cL{{\mathcal L}}
\def\cM{{\mathcal M}}
\def\cN{{\mathcal N}}
\def\cO{{\mathcal O}}
\def\cP{{\mathcal P}}
\def\cQ{{\mathcal Q}}
\def\cR{{\mathcal R}}
\def\cS{{\mathcal S}}
\def\cT{{\mathcal T}}
\def\cU{{\mathcal U}}
\def\cV{{\mathcal V}}
\def\cW{{\mathcal W}}
\def\cX{{\mathcal X}}
\def\cY{{\mathcal Y}}
\def\cZ{{\mathcal Z}}

\def \L{{\bbbl}}
\def \K{{\bbbk}}
\def \Z{{\bbbz}}
\def \N{{\bbbn}}
\def \Q{{\bbbq}}
\def\E{{\mathbf E}}
\def\G{{\mathcal G}}
\def\O{{\mathcal O}}
\def\cS{{\mathcal S}}
\def \R{{\bbbr}}
\def\\{\cr}
\def\({\left(}
\def\){\right)}

\def\Pb{\overline{P}}
\def\Qb{\overline{Q}}
\def\fl#1{\left\lfloor#1\right\rfloor}
\def\rf#1{\left\lceil#1\right\rceil}

\def \Prob{{\mathrm {}}}

\def \li {\mathrm {li}\,}

\newcommand{\comm}[1]{\marginpar{%
\vskip-\baselineskip 
\raggedright\footnotesize
\itshape\hrule\smallskip#1\par\smallskip\hrule}}

\title
[Products in  Residue Classes]{Products
in  Residue Classes}

\author[J. B.~Friedlander]{John B.~Friedlander}
\address{Department of Mathematics, University of Toronto,
         Toronto, Ontario M5S 2E4, Canada}
\email{frdlndr@math.toronto.edu}

\author[P. Kurlberg]{P\"ar Kurlberg}
\address{Department of Mathematics, Royal Institute of Technology,
SE-100 44 Stockholm, Sweden}
\email{kurlberg@math.kth.se}

\author[I. E.~Shparlinski]{Igor E.~Shparlinski}
\address{Department of Computing, Macquarie University, North Ryde,
Sydney, NSW 2109, Australia}
\email{igor@ics.mq.edu.au}

\begin{abstract}
We consider a problem of P.~Erd{\H o}s, A.~M.~Odlyzko and
A.~S{\'a}rk{\H o}zy about the representation of residue classes modulo $m$
by products of two not too large primes. While it seems that even the
Extended
Riemann Hypothesis is not powerful enough to achieve the expected
results, here we obtain some unconditional results ``on average'' over
moduli $m$ and residue classes modulo $m$ and somewhat stronger results
when the average is restricted to
prime moduli $m = p$. We also consider
the analogous question wherein the primes are replaced by
easier sequences so,
quite naturally, we obtain much stronger results.
\end{abstract}

\maketitle

\section{Introduction}

Let  $m$ and $a$ be integers with  $m\ge 1$ and $\gcd(a,m) =1$.
For $\cR$ and $\cS$, sets of positive integers, we consider the question
of whether there are integers $r\in \cR$ and $s\in \cS$ such that
$$
rs \equiv a \pmod m
$$
and, if so, how small can we choose these factors to be.

An obvious greedy algorithm would be to choose some small value
for one of these, say $r$,
and then look for the least $s\in \cS$ which satisfies the congruence
\begin{equation}\label{eq:greed}
s \equiv a\overline r \pmod m\ ,
\end{equation}
where $\overline r$ indicates the multiplicative inverse of $r$ modulo $m$.
It is clear that the use of this strategy  usually limits the possibility
for choosing $s$ to a range $s\le m^{1+ o(1)}$. Hence, to even obtain a
bound $r<m$, $s<m$ requires a more delicate argument and represents a
result of a different order of difficulty. Actually, one could even
hope to attain a bound
$$r, \, s\le m^{\frac12 + o(1)}\ ,
$$
a better result being hopeless, but such a goal seems far away
even in the simplest case where $\cR = \cS$ is the set of
all positive integers.

We are especially interested in the problem
where $\cR = \cS$ is the set of primes.
In this case we denote by $P(x;m,a)$
the number of solutions to the congruence
$$
p_1 p_2 \equiv a \pmod m
$$
in primes $p_1, p_2\le x$,
and we are interested in studying $P(x;m,a)$
for values of $x$, as small as possible, for example for $x = m$.

Our work has been motivated by the study of the quantity $P(x;m,a)$
in the paper~\cite{EOS} of P.~Erd{\H o}s, A.~M.~Odlyzko and
A.~S{\'a}rk{\H o}zy. They prove a number of results
conditional on various assumptions about the zero-free regions for
Dirichlet $L$--functions. However, even the Extended
Riemann Hypothesis\footnote{The assertion that  all nontrivial zeros of
  any Dirichlet $L$-function lie on the critical line.}
(ERH) seems not to be powerful enough
to prove their intended goal that
$$
P(m;m,a) > 0
$$
whenever $\gcd(a,m)=1$. Furthermore, even various relaxations
of this question considered in~\cite{EOS} have required
some unproven assumptions. In particular, under a certain weakened
form of the ERH, the mean-square
$$
P(x,m) = \sum_{\substack{a=1\\ \gcd(a,m)=1}}^m \(P(x;m,a) -
\frac{\pi(x)^2}{\varphi(m)}\)^2
$$
has been estimated successfully in the case $x=m$ and $m$ prime.
Here, as usual, $\varphi(k)$ is the Euler function and $\pi(x)$ is the number
of primes $p \le x$

We further relax the original question
but instead concentrate on unconditional results.
In particular, in Section~\ref{sec:Prime}
     we use the large sieve
to estimate $P(x,q)$ on average
over primes  $M<q\le 2M$ in ranges $M\ge x$
and also $P(x,m)$ over general integer moduli $m$.
In the case of the average over prime moduli we can come within a power of
a logarithm of the optimal range.

We also study the problem for some integer sets a little less
difficult than the primes. For
example, the sequence of squarefree integers is one which can be handled with
greater success and without any unproved assumptions. Let
$S(x;m,a)$
denote the counterpart of $P(x;m,a)$ wherein the primes $p_1, \,p_2$
are replaced by squarefree integers. Here, in
Section~\ref{sec:SF}, we obtain an asymptotic formula for
$S(x;m,a)$ which is nontrivial for $x \ge m^{3/4 + \varepsilon}$
for any fixed $\varepsilon > 0$ and sufficiently large $m$. As hinted
above, one might hope that such
formulas hold even down as far as $x \ge m^{1/2 + \varepsilon}$ but if
so this seems quite difficult. We
do not know how to get a wider range of uniformity (apart from the
$\varepsilon$) even for the
apparently easier problem where we do not insist that the
factors be squarefree!
In this case, where $\cR = \cS$ is the set of all positive integers,
the exponent $3/4$ rests on the Weil bound for Kloosterman sums
and has resisted improvement for half a century.
See however~\cite{Gar2,KhShp} for recent work related to this problem.

In Section~\ref{sec:P x SF} we consider the hybrid problem with
products $ps$ of a prime $p$ and a squarefree integer
$s$ in the range $p,s \le m^{1/2+\varepsilon}$
and show that, for any integer $m$,
these products  represent almost all reduced residue
classes modulo $m$ the expected number of times.

Finally, in Section~\ref{sec:Other} we consider a different example
wherein one of the two
factor sets is a sumset, a case in which rather general results can be
obtained provided neither set is very thin.
We also consider the case
of products of two primes and one shifted prime which
is also accessible
by present methods.

Throughout the paper the letters $p$ and $q$ are reserved for prime
numbers. The M\"obius function $\mu$ and
divisor function $\tau$ have their usual meanings.

\bigskip

\noindent{\bf Acknowledgements.}
Much of the work on this paper was done during visits by P.~K.\ and I.~S.\
to the University of Toronto, whose support and hospitality
are gratefully acknowledged.
Research of
J.~F.\ was partially supported by NSERC grant A5123, that of
P.K. by grants from the G\"oran
Gustafsson Foundation, and the Royal Swedish Academy of Sciences, and that
of I.~S.\ by ARC grant DP0556431.

\section{Products of Primes}
\label{sec:Prime}

We note that, for a
sufficiently large constant $c$, the result of Heath-Brown~\cite{HB2}
on the Linnik problem of the least prime in an arithmetic progression
implies by~\eqref{eq:greed}
that $P(x;m,a)>0$ for any $a$ such that $\gcd(a,m) =1$,
provided $x \ge cm^{11/2}$, and it has long been known that,
under the ERH, the exponent $11/2$ may be replaced by any number larger
than $2$.

It is even
expected that $x \ge m^{1 + \varepsilon}$ for any fixed $\varepsilon > 0$
is admissible but ideas for any
reasonable approach to this are lacking, at least for individual
progressions. However, one can show, again
using the greedy algorithm~\eqref{eq:greed} but now in conjunction
with the Barban-Davenport-Halberstam Theorem
(see~\cite[Theorem~17.2]{IwKow}),
that $P(x;m,a)>0$ for most $m\le M$ and most reduced
classes modulo $m$ provided that
$x/M(\log M)^3 \rightarrow \infty$.
However, we are able to do better than that with a different argument.

Let us define
$$
R(x,M) = \sum_{M<m \le 2M} P(x,m)\quad\text{and}\quad R_{\pi}(x,M)
= \sum_{M<q \le 2M} P(x,q)\ ,
$$
where, as usual, $q$ runs over primes.

We now improve on the trivial bounds:
$$
R(x,M)  \ll x^4 \quad\text{and}\quad R_{\pi}(x,M)  \ll x^4/\log x .
$$

\begin{theorem}
\label{thm:Prime Prod}
The following bounds hold:
$$
R(x ,M)\ll x^4(\log x)^{-A}+ Mx^2 ,
$$
for any $A$, with an implied constant that depends on $A$,
and
$$
R_{\pi}(x,M) \ll \(M^{-1}x^4 + Mx^2\)(\log x)^{-2}\ .
$$
\end{theorem}

\begin{proof}
Let $\cX_m$ be the set of all $\varphi (m)$ multiplicative characters
modulo $m$, and $\cX_m^*$
the set of primitive characters modulo $m$ (which in the case of
prime modulus $q$ includes all such
characters other than the principal character).

Using the orthogonality relation
$$
\frac{1}{\varphi (m)}
         \sum_{\chi \in \cX_m}\chi\(r\)
=\left\{\begin{array}{rll}1 &\text{if}\ r
\equiv 1 \pmod m, \\ 0 &\text{otherwise,}
\end{array}\right.
$$
for $\gcd(a,m)=1$, we write
\begin{eqnarray*}
P(x;m,a) &=& \sum_{p_1, p_2 \le x} \frac{1}{\varphi (m)}
         \sum_{\chi \in \cX_m}\chi\(p_1 p_2 a^{-1} \)\\
&=& \frac{\pi(x)^2}{\varphi (m)} +  \frac{1}{\varphi (m)}
\sum_{\substack{\chi \in \cX_m\\ \chi \neq \chi_0 }}  \chi\(a^{-1}\)
\sum_{p_1, p_2 \le x} \chi\(p_1 p_2 \)
\\
&=& \frac{\pi(x)^2}{\varphi (m)} +  \frac{1}{\varphi (m)}
\sum_{\substack{\chi \in \cX_m\\ \chi \neq \chi_0 }}
\chi\(a^{-1}\)  T_{\chi}(x)^2
\end{eqnarray*}
where
$$
        T_{\chi}(x) = \sum_{p  \le x} \chi\(p \).
$$
In particular,
\begin{eqnarray*}
P(x,m) &=&  \frac{1}{\varphi (m)^2}
\sum_{(a,m)=1}
\Biggl(\sum_{\substack{\chi \in \cX_m\\ \chi \neq \chi_0 }}
      \chi\(a^{-1}\) \Bigl(
\sum_{p  \le x} \chi\(p \)\Bigr)^2\Biggr)^2 \\
&\le&  \frac{1}{\varphi (m)^2}
\sum_{(a,m)=1}
\Biggl|\sum_{\substack{\chi \in \cX_m\\ \chi \neq \chi_0 }}
     \overline \chi\(a\)
        T_{\chi}(x)^2 \Biggr|^2\\
&=&  \frac{1}{\varphi (m)^2}
\sum_{(a,m)=1}
        \sum_{\chi_1, \chi_2 \neq \chi_0}\overline  \chi_1\(a\) \chi_2\(a\)
        T_{\chi_1}(x)^2  T_{\overline  \chi_2}(x)^2  \\
&=&  \frac{1}{\varphi (m)^2}
        \sum_{\chi_1, \chi_2 \neq \chi_0}
        T_{\chi_1}(x)^2  T_{\overline  \chi_2}(x) ^2
\sum_{(a,m)=1}\overline \chi_1\(a\)  \chi_2\(a\)\ ,
\end{eqnarray*}
where $\overline\chi$ denotes the conjugate character.
Since
$$
\sum_{(a,m)=1}\overline \chi_1\(a\) \chi_2\(a\)
=\left\{\begin{array}{rll} \varphi(m)  &\text{if}\  \chi_1 =\chi_2, \\
0 &\text{otherwise,}
\end{array}\right.
$$
we obtain
$$
P(x,m) \le
        \frac{1}{\varphi (m)}
        \sum_{\chi \neq \chi_0}
        T_{\chi }(x)^2  T_{\overline  \chi}(x) ^2
=  \frac{1}{\varphi (m)}
        \sum_{\chi \neq \chi_0}
        \left|T_{\chi }(x)\right|^4.
$$

We remark that
$$
T_{\chi}(x)^2 = \sum_{n \le x^2} a_n \chi(n)
$$
where $a_n = 2$ if $n =p_1p_2$ for two distinct primes
$p_1, p_2 \le x$, $a_n = 1$ if $n$ is the square of a prime $p\le x$,
     and $a_n = 0$ otherwise. Hence
$$
P(x,m) \le  \frac{1}{\varphi(m)}
        \sum_{\substack{\chi \in \cX_m\\ \chi \neq \chi_0 }}
        \Bigl|\sum_{n \le x^2} a_n \chi(n)\Bigr|^2
$$
which leads to the bound
\begin{equation}\label{eq:genmod}
        \sum_{M< m \le 2M} P(x,m)
\le  \sum_{M< m \le 2M}\frac{1}{\varphi(m)}
        \sum_{\substack{\chi \in \cX_m\\ \chi \neq \chi_0 }}
        \Bigl|\sum_{n \le x^2} a_n \chi(n)\Bigr|^2\ .
\end{equation}

We first treat the simpler
case where the average is over prime moduli $m=q$ so that all
non-principal characters modulo $q$ are primitive and~\eqref{eq:genmod}
can be replaced by the bound
$$
        \sum_{M < q \le 2M} P(x,q)
\ll \frac{1}{M} \sum_{M< q \le 2M}
        \sum_{\chi \in \cX_q^*}
        \Bigl|\sum_{n \le x^2} a_n \chi(n)\Bigr|^2\ .
$$
By the multiplicative form of the
large sieve inequality, see for example~\cite[Theorem~7.13]{IwKow},
we have
\begin{eqnarray*}
     \sum_{M< q \le 2M} \sum_{\chi \in \cX_q^*}
        \Bigl|\sum_{n \le x^2} a_n \chi(n)\Bigr|^2
& \ll & ( M^2 +  x^2) \sum_{n \le x^2} a_n^2 \\
& \ll & (M^2 +x^2)x^2 (\log x)^{-2}\ .
\end{eqnarray*}
Therefore,
$$
     R_{\pi}(x,M) =     \sum_{M < q \le 2M} P(x,q)
\ll  (Mx^2 + M^{-1}x^4)(\log x)^{-2}\ ,
$$
which concludes
the proof for the case  of prime moduli.

We now turn to the case of general modulus $m$ and need to
estimate, this time in general, the sum
$$S = \sum_{M< m \le 2M}\frac{1}{\varphi(m)}
        \sum_{\substack{\chi \in \cX_m\\ \chi \neq \chi_0 }}
        \Bigl|\sum_{n \le x^2} a_n \chi(n)\Bigr|^2
$$
on the right hand side of inequality~\eqref{eq:genmod}. Given a
character $\chi$ modulo $m$ occurring in this sum, let $\chi$ be
induced by a primitive
character $\psi$ modulo $f$ where $m=fe$ and, since $\chi$ is non-principal,
$f>1$. We have
$$\sum_{n \le x^2} a_n\chi(n)= \sum_{\substack{n \le x^2\\
\gcd(n,e)=1}} a_n \psi (n) =
\sum_{n \le x^2} a_n \psi (n)  +O(\log^{2} e)
$$
in view of the definition of $a_n$. Using this
and the inequality $\varphi(fe)\ge \varphi(f)\varphi(e)$, we have
\begin{equation*}
\begin{split}
S & \le \sum_{e\le 2M} \frac{1}{\varphi(e)}
\sum_{2\le f\le 2M/e} \frac{1}{\varphi(f)} \sum_{\psi \in \cX_f^*}
\Bigl|\sum_{n\le x^2} a_n \psi (n)\Bigl|^2 +O(Mx^2) \\
& = \sum_{e\le 2M} \frac{1}{\varphi(e)}
\{S_e(f\le F) + S_e(f > F)\} +O(Mx^2)\ ,
\end{split}
\end{equation*}
say, where $F =(\log x)^B$ for some large fixed $B$
and $S_e(f\le F)$ and $S_e(f > F)$ are the parts of the inner sums
taken over $f\le F$ and $f > F$, respectively.
(Note that we can
assume $\log M \ll \log x$ else the theorem is trivial.)
For $2\le f\le F$ we split the sum over $n$ into arithmetic progressions
modulo $f$ and apply to each of them the bound
$$\sum_{\substack{n\le x^2 \\ n\equiv b\, \pmod f}}a_n -
\frac{1}{\varphi(f)}\sum_{\substack{n\le x^2 \\ \gcd(n,f)=1}}a_n
\ll x^2(\log x)^{-C}
$$
for any $C$, which follows quickly from the Siegel-Walfisz theorem.
Using orthogonality the main term in the sum over $n$ disappears
and we obtain for each $\psi$ the bound
$$
\sum_{n\le x^2} a_n \psi (n)\ll  Fx^2 (\log x)^{-C}
$$
   from which we derive
$$
S_e(f\le F)\ll F^3x^4 (\log x)^{-2C}\ .
$$
For the sum over $f> F$ we split the sum into $\ll \log 2M$
dyadic intervals $(V,2V]$ and
apply again the same large sieve inequality as we did in the case
of prime moduli. We obtain
\begin{equation*}
\begin{split}
S_ e(f> F) & \ll \log 2M
\sup_{F\le V\le 2M}V^{-1}(V^2+x^2)\sum_{n \le x^2}|a_n|^2\\
& \ll (M +F^{-1}x^2) x^2  \ll M x^2  + x^4 (\log x)^{-B}\ .
\end{split}
\end{equation*}
We take $C=2B$,\, $B=A+1$, and sum over $e$, completing the proof.
\end{proof}

We remark that
the Siegel-Walfisz theorem restricts us to choose $F$ no larger
than a fixed power of $\log x$ which limits the saving in
     Theorem~\ref{thm:Prime Prod}
in this case of general modulus.

Let $W(M,x)$ be the number of pairs $(q,a)$ where the prime $q$
and integer $a$ satisfy
$M <q\le 2M$ and $1 \le a < q$ and such that $P(x;q,a)=0$.
Then
$$
     \frac{\pi(x)^4}{4M^2} W(M,x)  \le  R_{\pi}(x,M).
$$
Hence by Theorem~\ref{thm:Prime Prod} we have
$$
W(M,x) = o(M\pi(M))
$$
for any $x\le M$ satisfying $xM^{-1/2}(\log M)^{-3/2}\to \infty$, which
is thus within a power of the logarithm of being best possible.
This may be compared with a result of
M. Z. Garaev~\cite{Gar1} wherein a better power of the logarithm is obtained,
but for products of integers, not necessarily prime.

Finally, taking $x=M$ we see that,
$$
W(M,M)  \ll M (\log M)^2.
$$

\section{Products of Squarefree Integers}
\label{sec:SF}

As in the case of products of primes we can quickly deduce some
bound by appealing to the known results, in this case for
the smallest square-free integer in an arithmetic progression.
Thus, from the result of Heath-Brown~\cite{HB1} on that problem
it follows trivially that $S(x;m,a)>0$ for any $a$ such that
$\gcd(a,m)$ is square-free,
provided  $x \ge m^{13/9 + \varepsilon}$.

\begin{theorem}
\label{thm:SF Prod} For all integers $m\ge 1$ and $a$ with $\gcd(a,m) =1$
and real positive $x$, we have
$$
S(x;m,a) = \frac{36}{\pi^4} \cdot \frac{x^2}{m} \prod_{p\mid m} \(1 +
\frac{1}{p}-\frac{1}{p^2}+ \frac{1}{p^3}\)^{-1} + O\(x m^{-1/4 + o(1)}\),
$$
where the product is taken over all prime numbers $p \mid m$.
\end{theorem}

We remark that, as stated in the introduction, this gives the
asymptotic formula in the range $x \ge m^{3/4 + \varepsilon}$
for fixed positive $\varepsilon$. We take $x<m$ in the proof.
A very slight modification is needed for larger $x$.

\begin{proof}
For real $U$ and $V$ we demote by
$N(U,V;m,b)$ the number of solutions to the congruence
$uv\equiv b \pmod m $ in positive integers $u\le U$, $v \le V$.

Recall that $\mu(d)$ denotes the M\"obius function.
By the inclusion-exclusion principle, we write
\begin{equation}
\label{eq:S and N}
S(x;m,a) = \sum_{\substack{d,e=1\\ \gcd(de,m)=1}}^\infty \mu(d) \mu(e)
N(x/d^2,x/e^2;m,ad^{-2}e^{-2}) .
\end{equation}
A standard application of bounds for incomplete Kloosterman sums
(see~\cite[Corollary~11.12]{IwKow})
leads to the asymptotic formula
\begin{equation}
\label{eq:asymp}
N(U,V;m,b) = UV \frac{\varphi(m)}{m^2}
+ O(m^{1/2 +  o(1)})
\end{equation}
uniformly over integers $b$ with $\gcd(b,m) =1$,
see~\cite{BeKh,FuKi} and references therein.

For $\tau(w)$,
the number of positive divisors of $w$,
we recall the well known bound
\begin{equation}
\label{eq:tau}
\tau(w) = w^{o(1)},
\end{equation}
see for example~\cite[Section~I.5.2]{Ten}.

We define two quantities
\begin{equation}
\label{eq: y and z}
y = x m^{-3/4} \qquad \text{and}\qquad z = x m^{-1/2}\ ,
\end{equation}
which will feature in the proof.

We use the asymptotic formula~\eqref{eq:asymp} for $de \le y$.
which after substitution in~\eqref{eq:S and N} yields
\begin{eqnarray*}
\lefteqn{S(x;m,a) =   x^2 \frac{\varphi(m)}{m^2} \sum_{\substack{de\le y\\
\gcd(de,m)=1}}
\frac{ \mu(d) \mu(e)}{d^2 e^2} }\\
& & \qquad   +~O\Biggl( y m^{1/2 +  o(1)} +
\sum_{\substack{de > y\\ \gcd(de,m)=1}} N(x/d^2,x/e^2;m,ad^{-2}e^{-2}) \Biggr)
\end{eqnarray*}
(since by~\eqref{eq:tau} there are at most $y^{1 + o(1)} = y m^{o(1)}$
such pairs $(d,e)$ in the first sum).

Using~\eqref{eq:tau}, for the first term we obtain
\begin{eqnarray*}
        \sum_{\substack{de\le y\\
\gcd(de,m)=1}} \frac{ \mu(d) \mu(e)}{d^2 e^2}
& = &
\sum_{\substack{d,e=1\\ \gcd(de,m)=1}}^\infty \frac{ \mu(d) \mu(e)}{d^2 e^2}
+ O\(\sum_{k\ge y}\frac{ \tau(k) }{k^2}\)\\
&= &  \Biggl(\sum_{\substack{d=1\\ \gcd(d,m)=1}}^\infty
\frac{ \mu(d) }{d^2}\Biggr)^2  + O(y^{-1 + o(1)}).
\end{eqnarray*}

We have
\begin{equation*}
        \sum_{\substack{d=1\\ \gcd(d,m)=1}}^\infty
\frac{ \mu(d) }{d^2}  =
        \prod_{p\nmid m} \(1 - \frac{1}{p^2}\)
= \zeta(2)   \prod_{p\mid m} \(1 - \frac{1}{p^2}\)^{-1}\ ,
\end{equation*}
where $\zeta(s)$ is the Riemann zeta-function.
Therefore
\begin{equation}
\label{eq:prelim}
\begin{split}
S(x;m,a)
        =  \frac{36}{\pi^4}x^2 &\frac{\varphi(m)}{m^2} \prod_{p\mid m} \(1 -
\frac{1}{p^2}\)^{-1}\\ & + O\(y m^{1/2 +  o(1)} + x^2 y^{-1} m^{-1+o(1)}  +
\sum_{j = 0}^J \Delta_j\),
\end{split}
\end{equation}
where
$$
J = \rf{2 \log (x/z)}
$$
and

\begin{eqnarray*}
\Delta_0 &=& \sum_{\substack{z\ge  de > y\\ \gcd(de,m)=1}}
N(x/d^2,x/e^2;m,ad^{-2}e^{-2}),\\
\Delta_j & = & \sum_{\substack{2^j z\ge  de > 2^{j-1} z\\ \gcd(de,m)=1}}
N(x/d^2,x/e^2;m,ad^{-2}e^{-2}), \quad j = 1,  \ldots J.
\end{eqnarray*}

To estimate $\Delta_0$,  we note that if
$u v \equiv ad^{-2}e^{-2} \pmod m$ and with $1 \le u \le x/d^2$ and
$1 \le v \le x/e^2$
then for each fixed pair $(d,e)$,  the product
$w = uv \le x^2d^{-2}e^{-2}$ belongs to a
prescribed residue class modulo $m$ and thus takes at most
$ x^2d^{-2}e^{-2}m^{-1} + 1$ possible values. In turn, each value
of $w$ gives rise to $\tau(w) = w^{o(1)} = m^{o(1)}$ pairs
$(u,v)$ with $uv = w$, see~\eqref{eq:tau}.
Therefore
\begin{eqnarray*}
\Delta_0 &\le & \sum_{\substack{z\ge  de > y\\ \gcd(de,m)=1}}
\(\frac{x^2}{d^{2} e^{2} m} + 1\)  m^{o(1)}
= x^2m^{-1 + o(1)} \sum_{ de \ge y}  \frac{1}{d^2 e^2 }  + z m^{o(1)}\\
         &=&  x^2m^{-1 + o(1)}  \sum_{ k \ge y }  \frac{\tau(k)}{k^2}  + z
m^{o(1)} =  x^2y^{-1}m^{-1 +
o(1)} + z m^{o(1)}.
\end{eqnarray*}

To estimate $\Delta_j$, with $j\ge 1$,   we note that $\Delta_j$
does not exceed the number of pairs $(d,e)$ of positive
        integers such that $de \le 2^j z$, $\gcd(de, m) =1$ and
$ad^{-2}e^{-2} \equiv w \pmod m$ for some positive integer
$w\le W_j$ where
$$
W_j =   4 \frac{x^2}{2^{2j} z^2},   \qquad j = 1,  \ldots J .
$$
Furthermore,
due to our choice of $z$, see~\eqref{eq: y and z},
we have $W_j \le m$.
Thus if
$d$ and $e$ are fixed, then solutions to the congruence
$uv\equiv ad^{-2}e^{-2} \equiv w \pmod m $,
where $1 \le w < m$, in positive integers $u\le x/d^2$, $v \le
x/e^2$ satisfy the equation $uv = w$. Hence, by~\eqref{eq:tau},
every pair $d$ and $e$ leads to
$\tau(w) = m^{o(1)}$ possible pairs $(u,v)$.
Collecting together $d$ and $e$ with the same value of
$de=k$, and using~\eqref{eq:tau} again,
$$
\Delta_j
\le  m^{ o(1)} T (2^jz,W_j), \qquad j = 1, \ldots, J,
$$
where $T(K,W)$ is the number of positive
        integers $k \le K$, $\gcd(k, m) =1$ and
$ak^{-2} \equiv w \pmod m$ for some positive
integer $w \le W$.

Using exactly the same arguments as, for example
in~\cite[Lemma~1]{Shp1}
(Fourier expansion of the remainder term in the counting function, completion
of the relevant exponential sum, and finally application of Weil's theorem to
the completed sum), we obtain
that if $W < m$ then
$$
T(K,W) =  \frac{W}{m}
\sum_{\substack{k=1\\ \gcd(k,m) = 1}}^K  1+ O\(m^{1/2+o(1)}\)
\le \frac{KW}{m}  + O\(m^{1/2+o(1)}\),
$$
        which in turn implies that
$$
        \Delta_j \le  m^{o(1)}\( \frac{x^2}{2^{j} z m} + m^{1/2}\),\qquad j
= 1,  \ldots, J.
$$
Therefore
\begin{eqnarray*}
\sum_{j=1}^J  \Delta_j & \le & x^{2} z^{-1} m^{-1+ o(1)}\sum_{j=1}^J
\frac{1}{2^{j}}  +  J
m^{1/2+o(1)} \\
& \le &  x^{2} z^{-1} m^{-1+ o(1)} + m^{1/2+o(1)}.
\end{eqnarray*}
        Substituting the bounds on  $\Delta_0$ and on $\Delta_j$,  $j = 1,
\ldots, J$,
in~\eqref{eq:prelim} we obtain
\begin{equation*}
\begin{split}
S(x;m,a)
        =  \frac{36}{\pi^4}x^2 &\frac{\varphi(m)}{m^2} \prod_{p\mid m} \(1 -
\frac{1}{p^2}\)^{-1}\\ & + O\(y m^{1/2 +  o(1)} + x^2 y^{-1} m^{-1+o(1)}
+   z m^{o(1)}\).
\end{split}
\end{equation*}
Recalling the choice~\eqref{eq: y and z}
of $y$ and $z$, we conclude the proof.
\end{proof}

\section{Products of Primes and  Squarefree Integers}
\label{sec:P x SF}

Here we study  congruences
$ps \equiv a \pmod m$
in primes $p\le x$ and squarefree integers $s\le x$.

In fact our approach works for congruences
$rs \equiv a \pmod m$
where $r \le x$ is an element of a very
general set $\cR$ with $\gcd(r,m) =1$ and $s\le x$ is squarefree.
Accordingly, we write $Q(\cR,x;m,a)$ for the number
of solutions of such  a congruence.

\begin{theorem}
\label{thm:P x SF} For all
positive integers $m$, real $x$ and sets
    $\cR \subseteq[1, x]$   of integers
$r$  with $\gcd(r,m)=1$, we have
$$\sum_{\substack{a=1\\ \gcd(a,m) = 1}}^m
\left|Q(\cR,x;m,a)  -  \vartheta_m  \frac{ |\cR|x}{m}
\right|
\le
   |\cR|^{3/4} x^{3/4}  m^{1/4 + o(1)} ,
$$
where
$$
\vartheta_m = \frac{6}{\pi^2}
\prod_{p \mid m} \(1 -\frac{1}{p^2}\)^{-1}.
$$
\end{theorem}

\begin{proof} Since
$$
\sum_{\substack{a=1\\ \gcd(a,m) = 1}}^m
Q(\cR,x;m,a) \le |\cR| x\ ,
$$
we see that unless
\begin{equation}
\label{eq:triv}
|\cR|x > m
\end{equation}
the bound is trivial.

Let  $U(\cR,y;m,a)$ be the number of
solutions to the congruence
$ru \equiv a \pmod m$
in  $r \in \cR$ and positive integers $u\le y$.
Our main tool is the bound
\begin{equation}
\label{eq:Aver U}
\sum_{a=1}^m \left|U(\cR,y;m,a) -    \frac{|\cR|y}{m}  \right|^2
\le |\cR| x   m^{o(1)},
\end{equation}
for $y\le x$,
which has been given in~\cite{Shp2}.

For every positive integer $d\le x^{1/2}$ we denote by
$V_d(\cR,x;m,a)$ the number of
solutions to the congruence
$rv \equiv a \pmod m$
in  $r \in R$ and positive integers $v\le x$ with
$v \equiv 0 \pmod {d^2}$.

Using the inclusion-exclusion principle, we write
\begin{equation}
\label{eq:Q and V}
Q(\cR,x;m,a) = \sum_{d=1}^\infty \mu(d) V_d(\cR,x;m,a),
\end{equation}
where, as before, $\mu(d)$ is the M\"obius function.

Clearly, if $\gcd(a,m)=1$ but $\gcd(d,m)>1$, then $V_d(\cR,x;m,a) = 0$.
Furthermore, for  $\gcd(ad,m)=1$ we have
$$
    V_d(\cR,x;m,a) =  U(\cR,x_d;m,a_d),
$$
where $x_d = \fl{x/d^2}$ and  $a_d$   is defined by the congruence
$a_d d^2 \equiv a
\pmod m$, $1 \le a_d < m$.

We now choose some parameter $z\ge 1$, to be specified later
and write~\eqref{eq:Q and V}
\begin{eqnarray*}
    Q(\cR,x;m,a) & = &
\sum_{\substack{d=1\\ \gcd(d,m)=1}}^\infty \mu(d)U(\cR,x_d;m,a_d)\\
& = &
\sum_{\substack{d \le z\\ \gcd(d,m)=1}}  \mu(d)U(\cR,x_d;m,a_d)\\
& & \qquad\qquad     +~O\(  \sum_{x^{1/2} \ge d > z} U(\cR,x_d;m,a_d)\)\\
& = &
\sum_{\substack{d \le z\\ \gcd(d,m)=1}}  \mu(d)     \frac{ |\cR|x_d}{
m} +O\(\sigma_1(a)+
\sigma_2(a) \),
\end{eqnarray*}
where
\begin{eqnarray*}
\sigma_1(a) & = &  \sum_{\substack{d \le z\\ \gcd(d,m)=1}}
\left|U(\cR,x_d;m,a_d) -
\frac{|\cR|x_d}{m}
\right|\\
\sigma_2(a) & = &   \sum_{\substack{x^{1/2} \ge d >z\\
\gcd(d,m)=1}}    U(\cR,x_d;m,a_d).
\end{eqnarray*}
As in the proof of Theorem~\ref{thm:SF Prod}, for the main term we obtain
\begin{eqnarray*}
\sum_{\substack{d \le z\\ \gcd(d,m)=1}}  \mu(d)     \frac{ |\cR|x_d}{ m}
& = &  \frac{ |\cR|x}{ m}\sum_{\substack{d \le z\\ \gcd(d,m)=1}}
\mu(d)     \frac{1}{d^2}
+ O\(\frac{z|\cR|}{m}\)\\
& = & \vartheta_m  \frac{ |\cR|x}{m} +  O\( \frac{ |\cR|x}{zm} +
\frac{z|\cR|}{m}\).
\end{eqnarray*}
Accordingly we obtain
\begin{equation}
\label{eq:Q and Sigmas}
\sum_{\substack{a=1\\ \gcd(a,m) = 1}}^m
\left|Q(\cR,x;m,a)  -  \vartheta_m  \frac{ |\cR|x}{m}\right|
    \ll \frac{ |\cR|x}{z} +  z|\cR|  + \Sigma_1+\Sigma_2,
\end{equation}
where
$$
\Sigma_1 =  \sum_{\substack{a=1\\ \gcd(a,m) = 1}}^m \sigma_1(a)
\qquad \text{and}\qquad
\Sigma_2 = \sum_{\substack{a=1\\ \gcd(a,m) = 1}}^m \sigma_2(a).
$$

Changing the order of summation and noticing that
due to the condition
$\gcd(d,m)=1$,
as $a$ runs through all the
reduced residue classes modulo $m$ then so does $a_d$, we write
\begin{eqnarray*}
\Sigma_1 & = &   \sum_{\substack{d \le z\\ \gcd(d,m)=1}}
\sum_{\substack{a=1\\ \gcd(a,m) = 1}}^m  \left|U(\cR,x_d;m,a) -
\frac{|\cR|x_d}{m}
\right|, \\
\Sigma_2 & = &   \sum_{\substack{x^{1/2} \ge d >z\\ \gcd(d,m)=1}}
     \sum_{\substack{a=1\\ \gcd(a,m) = 1}}^m U(\cR,x_d;m,a).
\end{eqnarray*}
By the Cauchy inequality and the bound~\eqref{eq:Aver U}, we have
\begin{equation}
\label{eq:Sigma1}
\Sigma_1 \le   \sum_{\substack{d  \le z\\ \gcd(d,m)=1}}
|\cR|^{1/2} x^{1/2} m^{1/2 + o(1)}
\le z   |\cR|^{1/2} x^{1/2} m^{1/2 + o(1)}.
\end{equation}

Furthermore, it is clear that
$$
     \sum_{\substack{a=1\\ \gcd(a,m) = 1}}^m   U(\cR,x_d;m,a) \le |\cR|
x_d \le |\cR| x/d^2.
$$
Therefore
\begin{equation}
\label{eq:Sigma2}
\Sigma_2 \ll  |\cR| x  \sum_{\substack{x^{1/2} \ge d >z\\ \gcd(d,m)=1}}
\frac{1}{d^2}
\ll \frac{ |\cR|x}{z}  .
\end{equation}
Substituting the bounds~\eqref{eq:Sigma1} and~\eqref{eq:Sigma2}
in~\eqref{eq:Q and Sigmas}, we derive
\begin{equation*}
\begin{split}
\sum_{\substack{a=1\\ \gcd(a,m) = 1}}^m
\left|Q(\cR,x;m,a)  -  \vartheta_m  \frac{ |\cR|x}{m}
\right|& \\
      \ll \frac{ |\cR|x}{z} +    z|\cR|&  + z   |\cR|^{1/2} x^{1/2}
m^{1/2 + o(1)}.
\end{split}
\end{equation*}
Clearly $z|\cR| \le z   |\cR|^{1/2} x^{1/2} m^{1/2}$, thus the second
term in the
last inequality can be dropped. Now taking
$z =  |\cR|^{1/4} x^{1/4} m^{-1/4}  $ and
remarking that~\eqref{eq:triv} implies that  $z>1$,
we conclude the proof.
\end{proof}

We see that
each reduced residue class modulo $m$
which contains no integer of the form $rs$ with    $r \in \cR$ and
a squarefree integer $s\le x$ contributes a term of
order
$ |\cR|x/m$
to the sum
estimated in Theorem~\ref{thm:P x SF}.

If we take $x = m^{1/2+ \varepsilon}$ for some fixed
$\varepsilon$ and  $\cR$ to be the set of primes $p \le x$
with  $p\nmid m$, we see that  the number of reduced classes  modulo $m$
which are not of the form $ps$ with a  prime  $p\le x$ and
a squarefree integer $s\le x$,
is at most
\begin{eqnarray*}
    |\cR|^{3/4} x^{3/4}  m^{1/4 + o(1)} m (|\cR|x)^{-1}
& \le &   |\cR|^{-1/4} x^{-1/4}  m^{5/4 + o(1)} \\
& \le & m^{1 -\varepsilon/2 + o(1)}  \le m^{1 -\varepsilon/3}\ ,
\end{eqnarray*}
provided that $m$ is large enough.

\section{Some Other Products}
\label{sec:Other}

Our work in this section is motivated by another nice, albeit conditional,
result in~\cite{EOS}, wherein it has been shown that if $\varepsilon>0$ and
$q$ is large enough, then
     all invertible elements
modulo $q$
can be written as a product of three primes $p_1,p_2,p_3 < q$ under
the assumption
that all Dirichlet $L$-series $L(s,\chi)$ are
non-vanishing for ${\mathrm{Re}}(s)>1-(3+\varepsilon) \log \log q/\log q$ and
$|{\mathrm{Im}}(s)|\leq q$, for all
non-trivial characters $\chi$ modulo $q$.

If we alter the problem slightly and consider numbers that are
products of a prime and a sum of two primes, that is numbers of the form
$p_1(p_2+p_3)$, we can unconditionally show that it suffices to take
$p_1,p_2,p_3 < q^{1-\delta}$ for
some $\delta>0$ to obtain all residues modulo $q$.
We also consider another modification, namely the
representation of residue classes
by the product $p_1p_2(p_3+b)$ of two primes and a shifted prime.

\begin{theorem}
\label{thm:sum prod}
      Let $q$ be a prime.  Given  any invertible
      element $a$ in $\Z/q\Z$, there are
$\pi(X)^3/q + O(q^{3/16 + o(1)} X^{57/32})$
      primes $p_1, p_2,
      p_3 \leq X<q$ such that $p_1(p_2+p_3) \equiv a \pmod q$.
\end{theorem}

\begin{proof}
We have
\begin{eqnarray*}
\lefteqn{
    |\{ p_1, p_2, p_3 \leq X \ :\ p_1(p_2+p_3) \equiv a \pmod q \}|}\\
     & & \qquad\qquad = \frac{1}{q} \sum_{t=1}^q \sum_{p_1,p_2,p_3 \leq X}
\eq( t( p_2+p_3 - a\overline p_1)  )\\
& & \qquad\qquad =  \frac{\pi(X)^3}{q} +
\frac{1}{q} \sum_{t=1}^{q-1} \sum_{p_1 \leq X}
\eq( - t a\overline p_1  )\sum_{p_2,p_3 \leq X}
\eq( t( p_2+p_3 )  ) .
\end{eqnarray*}
The error term is bounded by
\begin{eqnarray*}
\lefteqn{
\frac{1}{q} \max_{t \not \equiv 0 \pmod q}
\left| \sum_{p_1 \leq X} \eq( - t a\overline p_1  )
\right| \cdot
\sum_{t=1}^{q-1}\left|\sum_{p_2,p_3 \leq X}
\eq( t( p_2+p_3 )  )\right|}\\
     & & \qquad\qquad =
\frac{1}{q}\max_{t \not \equiv 0 \pmod q}
\left| \sum_{p_1 \leq X} \eq( t \overline p_1  )
\right| \cdot\sum_{t=1}^{q-1}
\left| \sum_{p_2 \leq X}\eq( tp_2  )\right|^2.
\end{eqnarray*}
By~\cite[Theorem~1.1]{foumic98} we have
\begin{equation}
\label{eq:Exp p}
\max_{t \not \equiv 0 \pmod q}
\left| \sum_{p_1 \leq X} \eq( t \overline p_1  )
\right| \le q^{3/16 + o(1)} X^{25/32}.
\end{equation}
Also, by orthogonality we get
$$
\sum_{t=1}^{q-1}\left|\sum_{p_2 \leq X}\eq( tp_2  ) \right|^2
= q \pi(X)
$$
which concludes the proof.
\end{proof}

Clearly the asymptotic formula of Theorem~\ref{thm:sum prod}
is nontrivial if $X \ge q^{38/39+\delta}$ for an arbitrary
     $\delta>0$ and sufficiently large $q$.

We remark that the
the proof of Theorem~\ref{thm:sum prod} extends
immediately to yield the following more general result.

\begin{theorem}
\label{thm:sum prod gen}
Let $q$ be prime and $(a,q)=1$. Let $\cR$ be any set of
positive integers
$r$ with $(r,q)=1$ and $\cS$ any set of integers $1\le s<q$.
Then we have
$$ \Biggl|\sssum_{\substack{r\in \cR\, s_1,s_2\in \cS \\ r(s_1+s_2)
\equiv a \,(\bmod q) }}1
-\frac1q |\cR||\cS|^2  \Biggr| \le |\cS|\max_{h \not \equiv 0 \,(\bmod q)}
\Bigl|\sum_{r\in \cR}\eq (h\overline r)\Bigr|\ .
$$
\end{theorem}

For example, in the case that $\cR$ is the set of primes, then
$\cS$ can be replaced by an arbitrary set of residue classes modulo $q$
satisfying $ |\cS | \ge q^{1-\delta}$ for some $\delta$ determined
by the results of~\cite{bou05b} or~\cite{foumic98}.

\smallskip
We now recall the bound of
     A.~A.~Karatsuba~\cite{Kar} which asserts that
if $q$ is prime and $X \ge q^{1/2+ \varepsilon}$
for some $\varepsilon$, then for any fixed
integer $b$ with $\gcd(b,q)=1$,
\begin{equation}
\label{eq:Char p}
\max_{\chi \in \cX_q^*} \left| \sum_{p \le X}
\chi(p+b) \right| \ll X^{1-\delta},
\end{equation}
where, as before, $\cX_q^*$ is the set of all
nontrivial multiplicative characters modulo $q$
and $\delta> 0$ depends only on $\varepsilon$.

Using this bound in the
the same way as in~\cite{FrShp}, one easily gets the following
result:

\begin{theorem}
\label{thm:prod shift}
      Let $q$ be a prime and let $b$ be an integer with
$\gcd(b,q)=1$.  There are two absolute constants $\eta,\kappa>0$
such that for $q > X \ge q^{1-\eta}$,  for any invertible
      element $a$ in $\Z/q\Z$, there are
$(1 + O(q^{-\kappa}))\pi(X)^3/q$
      primes $p_1, p_2,
      p_3 \leq X<q$ such that $p_1p_2(p_3+b) \equiv a \pmod q$.
\end{theorem}

\section{Remarks}
We note that~\cite[Theorem~A.9]{bou05b} and~\cite[Corollarie~1.6]{foumic98}
give nontrivial estimates for the exponential sums in~\eqref{eq:Exp p}
as long as $ X > q^{1/2+\varepsilon}$  and
$ X>q^{3/4+\varepsilon}$, respectively. However these bounds are
less explicit than~\eqref{eq:Exp p} and thus only lead to a weaker
inexplicit statement.

In principle, Theorem~\ref{thm:prod shift} can
be extended
to composite moduli $m$ and more general products. In fact,
Z.~Kh.~Rakhmonov~\cite{Rakh} provides an analogue of the
bound~\eqref{eq:Char p},
however only for $X \ge m^{1 + \varepsilon}$.

We have already remarked on the work of M.~Z.~Garaev~\cite{Gar1}
showing that, in the case that $\cR = \cS$ is the set of all integers,
one can come within a small power of the logarithm of the expected
conjecture for most residue classes to most prime moduli. In very
recent work, M.~Z.~Garaev and A.~A.~Karatsuba~\cite{GarKar} have
proved that in fact this holds for {\it all} prime moduli
(see~\cite{GarGar,Shp2} for various refinements of this result).
Our bound for $R_{\pi}(x,M)$ in Theorem~\ref{thm:Prime Prod} gives
in particular an analogue of the original result~\cite{Gar1} in
the apparently harder case when $\cR = \cS$ is the set of primes.
However, obtaining an analogue to the result of~\cite{GarKar}
for the case of prime products remains an open problem.

\end{document}